\documentclass[12pt,leqno]{article}
\usepackage{amsmath,amssymb,latexsym}

\allowdisplaybreaks

\def\nek{,\ldots,}
\def\NN{\mathbb N}
\def\RR{\mathbb R}

\def\ZZ{\mathbb Z}
\def\alp{\alpha}
\def\bet{\beta}
\def\del{\delta}
\def\Del{\Delta}
\def\intl{\int\limits}
\def\prodl{\prod\limits}
\def\suml{\sum\limits}
\def\gam{\gamma}
\def\Gam{\Gamma}

\def\kap{\kappa}
\def\const{{\rm const}}

\def\Sig{\Sigma}
\def\calL{\cal L}
\def\done{{1\hskip-2.5pt{\rm l}}}
\def\qed{\Box}
\def\Fix{\hbox{Fix} (f)}

\newcommand{\Diffo}{{\text{Diff}_0}}
\newcommand{\p}{{\partial}}

\def\bks{{\backslash}}
\def\tet{\theta}

\def\nek{,\ldots,}

\newtheorem{theorem}{Theorem}[section]
\newtheorem{lemma}[theorem]{Lemma}
\newtheorem{cor}[theorem]{Corollary}
\newtheorem{question}[theorem]{Question}

\begin{document}
\title{A growth gap for diffeomorphisms of the interval}
\author{Leonid Polterovich and Mikhail Sodin\\
School of Mathematical Sciences\\
Tel Aviv University\\
Ramat Aviv, Israel  69978\\
polterov@post.tau.ac.il \quad sodin@post.tau.ac.il
}

\date{}
\maketitle

\begin{abstract}
Given an orientation-preserving diffeomorphism of the interval
$[0;1]$, consider the uniform norm of the differential of
its $n$-th iteration. We get a function of $n$ called the growth sequence.
Its asymptotic behaviour is an interesting invariant which
naturally appears both in geometry of the
diffeomorphisms groups and in smooth dynamics.
Our main result is the following Gap Theorem:
the growth rate of this sequence is either  exponential, or
at most quadratic with $n$. Further, we
construct diffeomorphisms with quite an irregular
behaviour of the growth sequence. This construction
easily extends to arbitrary  manifolds.
\end{abstract}

\baselineskip=18pt
%section 1
\section{Introduction and main results}\label{sec:intro}

Denote by $\Diffo ([0;1])$ the group of all $C^1$-diffeomorphisms $f$
of the interval $[0;1]$ such that $f(0)=0$ and $f(1)=1$. Given a
diffeomorphism $f\in \Diffo ([0;1])$, define its {\it growth
sequence\/}
$$
\Gam_n(f)=\exp ||\log (f^n)'||_\infty  = \max\left(
||(f^n)'||_\infty ,\, ||(f^{-n})'||_\infty \right)\,, \qquad n\in
\NN\ .
$$
Here $||\, . \,||_\infty$  stands for the uniform norm, and $f^n$,
$n\in \ZZ$, denotes the $n$~-th iterate of $f$. Let us say that
two sequences of positive real numbers are {\it equivalent\/} if
their ratio is bounded away from $0$ and $+\infty$.  The
equivalence class of the sequence $\Gam_n(f)$ is called the {\it
growth type\/} of $f$. Clearly, it is invariant under
conjugations in the group of diffeomorphisms.

From the viewpoint of dynamics, the growth type reflects
asymptotic distortion of length under iterations of $f$.
Geometrically, the growth type of $f$ is closely related to the
distortion of the cyclic subgroup $\{ f^n\}\subset \Diffo ([0;1])$
with respect to the multiplicative norm $\Gamma_1(f) $ on $\Diffo
([0;1])$. In \cite{DG} D'Ambra and Gromov suggested to study the
growth type for various classes of diffeomorphisms.

The growth sequence is always submultiplicative:
\begin{equation}
\label{eq1.1}
\Gamma_{n+m} (f) \le \Gamma_n(f) \Gamma_m(f)\,,
\end{equation}
therefore there always exists the limit
$$
\gamma (f) = \lim_{n\to\infty} \Gamma_n^{1/n}(f)\,.
$$
Let $\Fix$ be the set of fixed points of $f$. Using standard
arguments of ergodic theory, it is easy to check that
\begin{equation}
\label{eq1.2} \gamma (f) = 1 \qquad \hbox{if\ and\ only\ if\ }
\qquad f'(\xi)=1 \quad \hbox{for \ every} \ \xi\in \Fix
\end{equation}
(we bring details in the end of Section~2). Otherwise, $\gamma (f)
>1$, so $\Gamma_n(f)$ grows exponentially fast. Loosely speaking,
the exponent $\gamma (f)$ distinguishes between the parabolic and
hyperbolic behaviour of diffeomorphisms. Our main result
establishes the growth gap between the parabolic and hyperbolic
cases.

\setcounter{theorem}{2}
\begin{theorem}[Growth gap]
\label{theo1.3}
Let $f\in \Diffo ([0;1])$ be a $C^2$-diffeomorphism with $\gamma
(f)=1$. Then
$$
\Gamma_n(f) \le {\rm Const}\, n^2
$$
for every $n\in\NN$.
\end{theorem}

As the proof shows, the $C^2$-condition can be relaxed by assuming
that $\log f'$ has bounded variation on $[0;1]$.\ Probably it
cannot be relaxed much further. Under the assumptions of the theorem, the
estimate is sharp. For instance, if $f\in \Diffo ([0;1])$ is a
$C^\infty$-diffeomorphism such that $\Fix = \{0,\,1\}$,
$f'(0)=f'(1)=1$ but $f''(0)\ne 0$, $f''(1)\ne 0$, then one can
check that the growth type of $f$ is $n^2$.

\medskip The result can be considered in the following more general
context. Let $G$ be a group endowed with a multiplicative
(pseudo)-norm $\rho$, that is a function $\rho: G\to [0;+\infty)$
satisfying $\rho (\done) =1$, $\rho(f)=\rho(f^{-1})$, and
$\rho(fg)\le \rho(f) \rho(g)$. By a {\em growth gap} we mean a gap
in the possible growth types of sequences $\rho(g^n)$ for $g\in
G$. Existence of growth gaps is known for finite-dimensional Lie
groups. As a toy model, consider the group $GL(m, \RR)$ endowed
with the operator norm. For instance, when $m=2$ the possible
growth types are given by $e^{cn}n^q$ where $c\ge 0$,
$q\in\{0,1\}$. Other examples of growth gaps are given by certain
discrete groups endowed with the norm $e^{l(w)}$, where $l(w)$ is
the word length of an element $w$ with respect to a chosen set of
generators. See \cite{LMR} for the treatment of lattices. As far
as we know, Theorem~\ref{theo1.3} gives the first example of a
growth gap for an infinite-dimensional Lie group (though see
\cite{PSib} for some steps in this direction in the context of
Hofer's metric on groups of area-preserving diffeomorphisms).

\medskip

Our second result starts with another observation:
\setcounter{equation}{3}
\begin{equation}
\label{eq1.5} \sum_{n\ge 1} \frac{1}{\Gamma_n(f)} < \infty
\qquad \hbox{for \ every}\  f\in \Diffo ([0;1])\setminus \{\done\}\,.
\end{equation}
Indeed, take a point $x_0\in [0;1]\setminus \Fix$ and assume, for
example, that $f(x_0)>x_0$. Put $x_n=f^n x_0$, $\delta_n = x_{n+1}
- x_n$. Note that $[x_0;x_1]=f^{-n} [x_n;x_{n+1}]$, so
\begin{equation}
\label{eq1.6}
\Gamma_n(f) \geq \max\limits_x(f^{-n})'(x)\ge
\del_0/\del_n.
\end{equation}
Obviously, $\suml_{n\in\ZZ}\del_n\le 1$ and thus (\ref{eq1.5})
follows. In particular, we see that
\begin{equation}
\label{eq1.7} \Gamma_n(f) \ge \hbox{const}\, n
\end{equation}
for ``most'' indices $n\in \NN$. In many cases, (\ref{eq1.7})
holds for {\it all} $n\in \NN$, see a brief discussion below.
However, the next theorem shows that  there are
non-trivial $C^\infty$-diffeomorphisms with an arbitrary
slow growth of $\Gamma_n(f)$ along a rare subsequence of indices
$n$.

Denote by $\calL$  the set of all strictly increasing sequences
$\{ u(n)\}$, $n\in\NN$,  of positive real numbers with $u (n)\to
+\infty$ as $n\to +\infty$.

\setcounter{theorem}{6}
\begin{theorem}\label{theo1.8}
For every sequence $u\in\calL$  there exists a diffeomorphism $f$
in $\Diffo([0;1])\bks\{\done\}$ such that
$$
\liminf_{n\to\infty} \frac{\Gamma_n(f)}{u(n)} \le 1\,.
$$
\end{theorem}

The diffeomorphism we construct in Theorem~\ref{theo1.8} must {\em
oscillate} near the end points. Consider the function
$v(x)=f(x)-x$. Assume for a moment that $v$  is monotone near
$x=1$.  After appropriate choice of $x_0$  we can think that $v$
is non-increasing on $[x_0;1]$. Then the sequence $\del_n =
f^{n+1}x_0 - f^n x_0$, $n\ge 1$, is monotone as well. Thus
$\del_n\le\frac{1}{n}$ since $\sum \delta_n \le 1$. Therefore, by
(\ref{eq1.6}), $f$ satisfies \eqref{eq1.7} for all indices $n\in
\NN$. Let us say that a diffeomorphism $f$ is {\it flat\/} at the
end points if $f'(p)=1$, $f^{(i)}(p)=0$ for $p\in \{ 0;1\}$ and
all integers $i\ge 2$. If $f$ is not flat near 1, then the Taylor
expansion of $f$ at the point $x=1$ yields monotonicity of the
function $v(x)$ for $x$ sufficiently close to $1$, and therefore
at least linear growth of the sequence $\Gamma_n(f)$. Hence the
diffeomorphism from Theorem~\ref{theo1.8} must be a flat one. In
what follows, we will design an oscillating function $v$ which
forces $f$ to satisfy $\Gam_{n_i}(f)=o(n_i)$, $\{n_i\}\subset\NN$.
Of course these oscillations are rare and small since $\Sig\del_n$
converges.  An additional difficulty is that they have to be not
too steep since we wish $f$  to be $C^\infty$-smooth.  Let us
mention that flat diffeomorphisms of $[0;1]$  with an oscillating
$v$ were considered in a different context in the papers
\cite{Se}, \cite{K}.

Theorem~\ref{theo1.8} has a straightforward corollary pertaining
to diffeomorphisms of arbitrary compact manifolds $M$. Let $\Diffo
(M)$ be the group of all $C^1$-smooth diffeomorphisms
isotopic to the identity map $\done$.  Given a diffeomorphism
$f\in \Diffo (M)$, define as above its {\it growth sequence\/}
$$
\Gam_n(f)=\max(\max\limits_{x\in M}\| d_x f^n\|,\,
\max\limits_{x\in M}\|d_xf^{-n}\|)\,, \qquad n\in \NN\ .
$$
Here $\| d_xf\|$  stands for the operator norm of the differential
$d_xf$  calculated with respect to a Riemannian metric on $M$.

\begin{cor}
\label{theo1.9} Let $B$  be a closed Euclidean ball of dimension
$\geq 1$. For every $u\in\calL$ there exists a $C^\infty$-diffeomorphism
$g\in \Diffo(B)\bks \{\done\}$ which equals the identity near the
boundary and whose growth sequence satisfies
\setcounter{equation}{8}
\begin{equation}
\label{eq1.11}
\liminf\limits_{n\to +\infty}\,
\frac{\Gam_n(g)}{u(n)}\le 1 \qquad \text{and} \qquad
\sum_{n=1}^{\infty} \frac{1}{\Gamma_n(g)} < \infty\;.
\end{equation}
\end{cor}

\medskip\par\noindent{\bf Proof: } Let $f$ be a diffeomorphism of
$[\frac1{3};\frac2{3}]$ which is flat at the end points and
satisfies $\liminf\limits_{n\to +\infty}\,
\frac{\Gam_n(f)}{u(n)}\le 1$. Existence of such diffeomorphism
follows from Theorem~\ref{theo1.8} and discussion above. Extend
$f$ by the identity map to the whole interval $[0;1]$. We get a
smooth diffeomorphism $h$ of $[0;1]$ which satisfies
(\ref{eq1.11}). Define now a diffeomorphism $g$ of the ball
$B=\{|x|\le 1\}$ by $g(x)=xh(|x|)/|x|$. Clearly, $g$ equals the
identity outside the spherical annulus $A=\{\frac1{3}\le |x|\le
\frac2{3} \}$. We claim that $g$ also satisfies conditions
(\ref{eq1.11}). Indeed, $g^n(x)=xh^n(|x|)/|x|$ for every
$n\in\ZZ$. Take a tangent vector $v\in T_x\RR^m$ and decompose it
as $v=\xi+\eta$, where $\xi$ is parallel to $x$ and $\eta$ is
orthogonal to $x$. One readily calculates that $$d_xg(v) =
(h^n)'(|x|)\xi + \frac{h^n(|x|)}{|x|}\eta$$ for $x\in B\setminus
\{0\}$. Note that $|x|^{-1}h^n(|x|)\in [\frac1{2}; 2]$ for all
$x\in A$. This yields
$$
\Gamma_n(h) \le \Gamma_n(g) \le \max(2, \Gamma_n(h)).
$$
Since $\Gamma_n(h)\to\infty$ in view of (\ref{eq1.5}), we see that
$\Gamma_n(g) = \Gamma_n(h)$ for large $n$. Hence $g$ also satisfies
conditions (\ref{eq1.11}). \hfill $\qed$

\medskip\par\noindent{\bf Remark. }
Corollary~\ref{theo1.9} immediately extends to an arbitrary
compact manifold. Indeed, fix a closed ball inside the manifold
and extend the diffeomorphism $g$  constructed in the theorem by
the identity map. We get a diffeomorphism from $\Diffo
(M)\setminus\{\done\}$ which satisfies (\ref{eq1.11}).

\smallskip
It is interesting to compare this remark with
restrictions on the growth type of symplectic maps which were
obtained recently in \cite{P1}, \cite{PS}. For instance, let $f$ be
an area-preserving $C^\infty$-diffeomorphism of the 2-dimensional torus
which
is isotopic to the identity map $\done$. Assume that $f\not=\done$
and $f$  has a fixed point.  Then according to \cite{PS} the
growth type of $f$ is at least linear, that is, (\ref{eq1.7})
holds for all indices $n\in \NN$. We refer to \cite{P1} for
extensions to other symplectic manifolds including all closed
surfaces of higher genus\footnote {P. LeCalvez informed us that he can
prove this result for surfaces by a different method. } and for further
discussion. Clearly the fixed point condition is crucial there.
Indeed, if $f$ is a translation of the torus, the sequence
$\Gam_n(f)$  is bounded (see \cite{P2} for more sophisticated
examples). On the other hand, as we have seen above, there exists
a $C^\infty$-diffeomorphism which has fixed points but violates inequality
\eqref{eq1.7} for a subsequence.

\subsection*{An Outlook}
%%%%%%%%%%%%%%%%%%%%%%%%%%%%%%%%%%%%%%%%%%%%%%%%
Flat fixed points form a major difficulty in the study of the
growth for $C^{\infty}$-diffeomorphisms of the interval. One can
speculate that further understanding of their influence on the
growth sequence leads to a rather satisfactory description of the
``spectrum" of possible growth types. The Gap Theorem and the
examples provided by Theorem \ref{theo1.8} correspond to the
opposite ends of this spectrum. To be more precise recall that if
a diffeomorphism $f$ has a non-degenerate fixed point $\xi$ (that
is $f'(\xi) \neq 1$) its growth is exponential. Assume now that
all fixed points of $f$ are degenerate but non-flat. We say that
$\xi \in \Fix $ has order $p \in \NN$ if $f^{(j)}(\xi) =0$ for all
$j=2,...,p$ but $f^{(p+1)}(\xi) \neq 0$. In this case one should
be able to show (using e.g. the Takens normal form \cite{T}) that
\begin{equation}
\label{eq1.10}
\Gamma_n(f) \sim n^{\frac{p+1}{p}},
\end{equation}
where $p$ is the minimal order of the fixed points of $f$.
Therefore, in the general case, we arrive at the following problem:
{\it What is the contribution of flat fixed points to the growth type of
$f$?} Warning: setting $p = \infty$ in formula
(\ref{eq1.10}) leads to the answer $ \Gamma_n \sim n$ which is obviously
wrong: it contradicts (\ref{eq1.5}).
In fact, Theorem \ref{theo1.8} shows that flat fixed points
sometimes yield an irregular behavior of the growth sequence.
Nevertheless an optimistic scenario is that
{\it the contribution of flat fixed points
does not exceed $n^{1+\epsilon}$
for every $\epsilon > 0$.} Note that the Gap Theorem
confirms this for $\epsilon = 1$.
If this is indeed true, we get
an infinite sequence of new gaps formed by the growth types
$$\left\{n^{\frac{p+1}{p}}\right\}, \qquad p \in \NN. $$
This picture, though quite an enticing , at the moment seems to be out of reach.
Our proof of the Gap Theorem goes in another direction and completely
ignores higher derivatives at the fixed points (as a compensation, we work
in the $C^2$-category). Let us conclude this discussion with the following
test
\setcounter{theorem}{10}
\begin{question}
Suppose that $f$ is a sufficiently smooth diffeomorphism of
$[0;1]$ with $\Fix = \{0;1\}$. Assume that $f'(0)=f'(1) = 1$ and
$f''(0) = f''(1)= 0$. Is it true that $\Gamma_n(f) = o(n^2)$ as
$n\to\infty$?
\end{question}

\smallskip\par\noindent\underline{Added on March, 2003.} In a
recent preprint \cite{B}, A.~Borichev confirms formula (1.10) and
gives the affirmative answer to Question~1.11. At the same time,
according to \cite{B} our ``optimistic scenario'' appeared to be
wrong.

%%%%%%%%%%%%%%%%%%%%%%%%%%%%%%%%%%%%%%%%%%%%%%%%

%section2
\section{Existence of the growth gap}
\setcounter{theorem}{0}
\setcounter{equation}{0}

In this section we prove Theorem~\ref{theo1.3} and check relation
(\ref{eq1.2}).

Let $f$ be a $C^2$-diffeomorphism of $[0;1]$ with $\gamma (f) =
1$. Let $v(f)$ be the variation of $\log f'$ on the interval
$[0;1]$. We shall use a classical

\begin{lemma}[Denjoy]
\label{lem2.1}
If $J\subset [0;1]$ is a closed interval such that $fJ\cap J = \emptyset$,
then for every $n\in \NN$ and every $x,y\in J$
$$
e^{-v(f)} \le \frac{(f^n)'(x)}{(f^n)'(y)} \le e^{v(f)}\,.
$$
\end{lemma}

\medskip
For $n\ge 0$ put
$$
a_n(f) = \max_{[0;1]} \log (f^n)'(x) \qquad \hbox{and} \qquad
a_n(f^{-1}) = \max_{[0;1]} \log (f^{-n})'(x)\,.
$$
Note that $a_0(f^\pm)=0$. These two sequences appear to be ``almost
convex'':

\begin{lemma}
\label{lem2.2}
The sequences $a_n=a_n(f)$ (resp. $a_n=a_n(f^{-1})$) satisfy
the inequality
$$
2a_n - a_{n-1} - a_{n+1} \le C(f)e^{-a_n}\,,
\qquad n\in \NN\,,
$$
with $C(f) = L(f)e^{v(f)}$, where $L(f)$ is the Lipschitz constant of
the function $\log f'$.
\end{lemma}

\medskip\par\noindent{\bf Proof of Lemma~\ref{lem2.1}:} We prove the
statement for the sequence $a_n=a_n(f)$. The proof for the second sequence
is the same (note that $v(f)=v(f^{-1})$).
Choose $x_0$ such that $a_n = \log (f^n)'(x_0)$. In view of (\ref{eq1.2})
$x_0$ does not belong to $\Fix$. Put $x_j=f^jx_0$, $j\in \ZZ$.  Then we
have
$$
a_{n+1} \ge \log (f^{n+1})'(x_{-1})
= \sum_{j=-1}^{n-1} \log f'(x_j)\,,
$$
$$
a_{n-1} \ge \log (f^{n-1})'(x_{1})
= \sum_{j=1}^{n-1} \log f'(x_j)\,,
$$
and
$$
a_{n} = \sum_{j=0}^{n-1} \log f'(x_j)\,.
$$
Therefore,
\begin{eqnarray*}
2a_n - a_{n-1} - a_{n+1} &\le& \log f'(x_0) - \log f'(x_{-1}) \\ \\
&\le& L(f) |x_0 - x_{-1}| \, < \, L(f)
\frac{|x_0-x_{-1}|}{|x_n-x_{n-1}|}
\\ \\
&=& \frac{L(f)}{(f^n)'(y)} \qquad \qquad \qquad \qquad \qquad
(y\in (x_{-1}, x_0)\,)
\\ \\
&\le& \frac{L(f)e^{v(f)}}{(f^n)'(x_0)} \, = \, C(f)e^{-a_n}\,,
\end{eqnarray*}
In the last inequality, we apply Lemma~\ref{lem2.1} to the interval
$[y;x_0]$. We are done. \hfill $\qed$

\medskip The next lemma is crucial:

\begin{lemma}[Growth lemma]
\label{lem2.3}
Let $\{a_n\}_{n \ge 1}$ be a sequence
of real numbers such that for each $n\ge 1$
\setcounter{equation}{3}
\begin{equation}
2a_n - a_{n-1} - a_{n+1} \le C\, e^{-a_n}\,,
\qquad C>0\,,
\label{eq1}
\end{equation}
and $a_0=0$. Then either for each $n\ge 1$
\begin{equation}
a_n \le 2\log \left(n\sqrt{\frac{C}{2}} +1\right)\,,
\label{eq2}
\end{equation}
or
\begin{equation}
\liminf_{n\to\infty} \frac{a_n}{n} > 0\,.
\label{eq3}
\end{equation}
\end{lemma}

\medskip\par\noindent{\bf Proof of Theorem~\ref{theo1.3}: } Since
$\Gamma_n(f) = \exp \left( \max (a_n(f), a_n(f^{-1}) \right)$,
Lemmas \ref{lem2.2} and \ref{lem2.3} yield
Theorem~\ref{theo1.3}. \hfill $\qed$

\medskip\par\noindent{\bf Proof of Lemma~\ref{lem2.3}:}
Introduce the second difference operator
$$L_j p = 2p_j - p_{j-1} - p_{j+1}\,, \qquad j \geq1\,,$$
acting on sequences $\{p_j\}$, $j \geq 0$.
Set $D=\sqrt{\frac{C}{2}}$ and observe that the sequence
$h_j=2\log(j\sqrt{C/2}+1)$ is a super-solution of the
non-linear second order difference equation
\begin{equation}
 \label{eqstar}
L_j p = Ce^{-p_j}\;.
\end{equation}
Indeed
$$L_j h = 2\log\frac{(Dj+1)^2}{(Dj+1)^2-D^2} >
2\frac{D^2}{(Dj+1)^2}= Ce^{-h_j}.$$
Here we used inequality $\log(A/B)> (A-B)/A$ valid
for $A > B > 0$. On the other hand condition (\ref{eq1})
says that the sequence $\{a_j\}$ is a sub-solution of the same equation
(\ref{eqstar}) with $a_0 = h_0 = 0$.
Our first claim is that for any $\epsilon \geq 0$ the sequence
$b_j = a_j - (1+\epsilon)h_j$ has no positive local maxima. This is
a version of the maximum principle for equation (\ref{eqstar}).
Indeed, suppose that $i \geq 1$ is a local maximum of $\{b_j\}$.
Then $L_i b \geq 0$, and
$$Ce^{-a_i} \geq L_i a = L_i b + (1+\epsilon)L_ih \geq L_i h \geq Ce^{-h_i}.$$
Hence $a_i \leq h_i$ and so $b_i \leq 0$. The claim follows.

Introduce the difference operator $\p _j p = p_j - p_{j-1}$, $j \geq 1$.
We shall show that if (\ref{eq2}) fails then
$\liminf\limits_{j \to+\infty} \p _j a$
is strictly positive, which clearly yields (\ref{eq3}).

If (\ref{eq2}) fails, there exists $m \in \NN$ and $\epsilon > 0$
so that $a_m > (1+\epsilon)h_m$. Fix $\epsilon$ and assume that $m$ is the
minimal positive integer which satisfies this inequality. This means that
$a_j \leq (1+\epsilon)h_j$ for $0 \leq j \leq m-1$.

Consider again the sequence $b_j = a_j - (1+\epsilon)h_j$.
By our assumption $b_m > 0$, and $b_j \leq 0$ for $0 \leq j \leq m-1$.
Since, as we checked above, this sequence cannot have positive local
maxima, we get that $b_n \geq b_{n-1}$ for all $n > m$.
Take any $n > m$. Then
\begin{equation}
\label{eq2.7}
\p_n a = \p_n b + (1+\epsilon) \p_n h > \p_n h \;.
\end{equation}
Furthermore, since
$$ a_n = b_n + (1+\epsilon)h_n > (1+\epsilon)h_n$$
we get
\begin{equation}
\label{eq2.9}
\p_n a - \p_{n+1}a = L_n a \leq Ce^{-a_n} \leq Ce^{-(1+\epsilon)h_n}\;.
\end{equation}
Take $N > n$ and observe that in view of (\ref{eq2.7}) and (\ref{eq2.9})
$$
\p_n h < \p_n a \leq Ce^{-(1+\epsilon)h_n} + \p _{n+1}a \leq \,...\,
\leq \sum_{j=n}^{N-1} Ce^{-(1+\epsilon)h_j} + \p_N a\,.
$$
Since the sequence $\{h_j\}$ is increasing
the first term on the right hand side does not exceed
$$
e^{-\epsilon h_n} \sum_{j=n}^{N-1} Ce^{-h_j} \leq  e^{-\epsilon h_n}
\sum_{j= n}^{N-1}
L_j h = e^{-\epsilon h_n} (\p _n h - \p_ N h) \leq e^{-\epsilon h_n} \p_n
h\,.
$$
This yields $\p _N a \geq (1 - e^{-\epsilon h_n} )\p_n h$.
Fix $n$ large enough so that $1 - e^{-\epsilon h_n} \geq 1/2$.
Letting $N\to +\infty$, we obtain that
$$\liminf\limits_{N \to +\infty} \p_N a \geq \frac{1}{2}\p_n h > 0 \;,$$
which yields (\ref{eq3}).
\hfill $\qed$

\medskip We complete this section with
\subsection*{Proof of relation (\ref{eq1.2})}
Suppose that $f'(\xi) = 1$ for each $\xi\in \Fix$. We have to show
that $\gamma (f) =1$ which means that $$\lim_{n\to\infty} a_n(f)/n
= \lim_{n\to\infty} a_n(f^{-1})/n = 0\,.$$ Put $F(x)=\log f'(x)$,
and assume on the contrary that
$$
\lim_{n\to\infty} \max_{x\in [0;1]} \left\{ \frac{1}{n}
\sum_{i=0}^{n-1} F(f^i(x))\right\} = \lim_{n\to\infty}
\frac{a_n}{n} = c>0\,.
$$
Using the Krylov-Bogolyubov argument, we choose a large enough $N$
and a point $y_{(N)}$ such that
$$
\frac{1}{N} \sum_{i=0}^{N-1} F(f^i y_{(N)} \,) > \frac{c}{2}\,,
$$
and consider a sequence of probability measures on $[0;1]$
$$
\sigma_N = \frac{1}{N} \sum_{i=0}^{N-1} \delta_{f^i y_{(N)}}\,,
$$
where $\delta_x$ is the Dirac measure at $x$. Then there is a
subsequence $N_j\to\infty$ such that $\sigma_{N_j}\to
\sigma_\infty$ (in the weak-* topology), and $\sigma_\infty$ is an
invariant measure of $f$ such that
$$
\int F\, d\sigma_\infty = \lim_{N_j\to\infty} \int F\,
d\sigma_{N_j} >0.
$$
Note that for every interval $J\subset [0;1]\setminus \Fix$ there
is a $k_0$ such that $f^kJ\cap J = \emptyset$ for every $k$ with
$|k|\ge k_0$. Hence the support of every invariant measure
$\sigma$ of $f$ must be contained in the set $\Fix$. Thus $\int
F\, d\sigma_\infty = 0$ since $F$ vanishes on $\Fix$ due to our
assumption. This contradiction proves that $\gamma(f)=1$. \hfill
$\qed$

%section3
\section{Diffeomorphisms with irregular growth sequences}\label{sec:cons}
\setcounter{theorem}{0}

In this section we prove Theorem~\ref{theo1.8}.
Fix a sequence $\{ u(n)\}\in\calL$  of positive real
numbers, $u(n)\nearrow +\infty$  as $n\to +\infty$.  For
a $C^\infty$-function $\Del:\RR\to (0;+\infty)$  define
recursively a sequence of functions
$$
g_0(t)=\frac{\Del(t+1)}{\Del(t)}\,,
\qquad
g_{m+1}(t)=\frac{g'_m(t)}{\Del(t)}\,.
$$

\setcounter{equation}{1}
\begin{theorem}\label{theo2.1}
There exists an even $C^\infty$-function $\Del:\RR\to
(0;+\infty)$  with the following properties:
%(2.2)
\begin{equation}\label{eq2.2}
\intl^{+\infty}_{-\infty}\Del(t)dt <\infty\ ;
\end{equation}
%(2.3)
\begin{equation}\label{eq2.3}
\text{there is a sequence of positive integers}\
\tau_i\nearrow +\infty\ \text{such that}
\end{equation}
$$\sup\limits_{t\in\RR}\frac{\Del(t+\tau_i)}{\Del(t)}\le
u(\tau_i)\ ;
$$
%(2.4)
\begin{equation}\label{eq2.4}
g_0(t)\to 1\ \text{and}\ g_m(t)\to 0\ \text{as}\
t\to\infty\ \text{for all integers}\ m\ge 1\ .
\end{equation}
\end{theorem}

First, we deduce Theorem~\ref{theo1.8}. Without loss of generality
assume that $\intl^{+\infty}_{-\infty}\Del(t) dt=1$. Put
$a(\eta)=\intl^\eta_{-\infty}\Del(s)ds$, and define
$f:[0;1]\to[0;1]$ by
$$
f(x)=\left[\begin{aligned}
&0,&\ \ x=0\phantom{(0;1)}\nonumber\\
&1,&\ \ x=1\phantom{(0;1)}\nonumber\\
&a(a^{-1}(x)+1),&\ \ x\in (0;1)\ .\nonumber
\end{aligned}
\right.
$$

\medskip\par\noindent{\bf Proof of Theorem \ref{theo1.8}: }
Let us verify that $f$  satisfies all the requirements of Theorem
\ref{theo1.8}.

1) We claim that $f$  is a $C^\infty$-diffeomorphism of
the closed interval $[0;1]$, and moreover $f$  is flat at the
end points.  Indeed, $f$
is continuous on $[0;1]$ and smooth on $(0;1)$.  Thus it
suffices to check that $f'(x)\to 1$  and $f^{(m)} (x)\to
0$, $m\ge 2$, as $x\to 0$  and $x\to 1$.
Put $\eta =a^{-1}(x)$.  Then $f'(x)=g_0(\eta)$,
$$f^{(2)}(x)=\frac{g_0'(\eta)}{\Del(\eta)}
= g_1(\eta)\nek
f^{(m+1)}(x)=\frac{g_{m-1}'(\eta)}{\Del(\eta)} = g_m(\eta)$$
for all $m\ge 2$.  Here we use that
$\frac{d\eta}{dx}=\frac{1}{\Del(\eta)}$.
The claim
follows from Property~\eqref{eq2.4} of $\Del$.

2) Note that $f^n(x)=a(a^{-1}(x)+n)$  for all $n\in\ZZ$.
Hence
$$\Gam_n(f)=\max \left(\max\limits_x\frac{a'(a^{-1}(x)+n)}
{a'(a^{-1}(x))}\ ,\ \max\limits_x\frac{a'(a^{-1}(x)-n)}
{a'(a^{-1}(x))}\right)= \sup\limits_{\eta\in\RR}\frac{\Del(\eta
+n)}{\Del(\eta)}$$ since $\Del$  is even.  Property \eqref{eq2.3}
guarantees that $\Gam_{\tau_i}(f)\le u(\tau_i)$. This completes
the proof.\hfill$\qed$

\medskip

It remains to prove Theorem~\ref{theo2.1}, that is to construct an
even $C^\infty$-function with properties
\eqref{eq2.2}--\eqref{eq2.4}.

\medskip\par\noindent{\bf Idea of the construction:}
As the first approximation for $\Delta$ we start with an even
$C^\infty$-function $h:\RR\to (0;+\infty)$
satisfying conditions (\ref{eq2.2}) and (\ref{eq2.4}) and
such that $h(0)=1$, and $h(t)$  decreases for $t>0$.
Consider the weighted average
$$
\left( A_{\tau, \mu} h\right)(t)
= \sum_{j\in \ZZ} \mu^{|j|} h(t+j\tau)\,,
$$
where $0<\mu<1$. It is not difficult to check that
$$
\mu \le \sup_{t\in \RR}
\frac{\left(A_{\tau,\mu} h\right)(t+\tau)}
{\left(A_{\tau,\mu} h\right)(t)} \le \frac{1}{\mu}\,,
$$
since the average produces ``self-similar humps'' of relative
magnitude $\mu$. Then iterating this procedure with properly
chosen sequences $\tau_i\to\infty$ and $\mu_i\to 0$, we get an
even function satisfying conditions (\ref{eq2.2}) and
(\ref{eq2.3}). Unfortunately, we loose the smoothness
property~(\ref{eq2.4}).

To mend this, we modify the operator $A$ by introducing an additional
rescaling:
$$
\sum_{j\in \ZZ}\mu^{|j|} h\big(\alpha_j (t+j\tau) \big)\,,
$$
where $\alpha_j$ are suitably chosen rescaling factors. The new average
still produces self-similar humps, this time without spoiling the
behaviour of
the derivatives. Then an infinite repetition of this process (with a
careful choice of values of $\tau$, $\mu$ and $\{\alpha_j\}$ in each step)
does the job.

\medskip\par\noindent{\bf Formal construction:} Let $h$ be an even
$C^\infty$-function $h:\RR\to (0;+\infty)$
such that $h(0)=1$, $h(t)$  decreases for $t>0$ and
$h(t)=\frac{1}{t\log^2t}$ for $t\ge 3$.
Take a function $\tau:\NN\to\NN$ such that
%(2.5)
\begin{equation}\label{eq2.5}
\suml_{i\in\NN}\frac{1}{\log u(\tau_i)} <\infty\ ;
\end{equation}
We also assume that the value $\tau_1$  is sufficiently
large.  Define now two functions $\mu:\NN\to (0;1)$  and
$\gam:\NN\times\ZZ\to (1;+\infty)$  as follows:
\begin{eqnarray*}
&&\mu_i=u^{-1/4}(\tau_i)\ ,\\
&&\gam_{i,\ell}=\min(|\log\mu_i|,\mu_i^{-\frac{1}
{\sqrt{|\ell|}}})\quad \text{for} \quad\ell\not= 0\,,
\\&&\text{and}\ \gam_{i,0}=|\log\mu_i|\ .
\end{eqnarray*}
Let $\ZZ^\infty$ be the space of all functions
$k:\NN\to\ZZ$  with finitely many non-zero values $k_i$.
Define two functions on $\ZZ^\infty$ as follows:
$$
\varphi(k)=\prodl_{i\in\NN}\mu^{|k_i|}_i\,, \qquad \tet(k)
=\prodl_{i\in\NN}\gam^{|k_i|}_{i,k_i}
$$
(the products have only finitely many factors).  Mention that the
functions $\varphi$  and $\varphi\tet$ are bounded by one:
$$
\varphi(k) \le \varphi(k)\tet(k) \le \prodl_{i\in \NN}
\mu_i^{|k_i|-\sqrt{|k_i|}} \le 1.
$$

\medskip\par\noindent
{\bf Fundamental definition:}  Put
$$\Del(t)=\suml_{k\in\ZZ^\infty}\varphi(k)h\big(\varphi(k)
\tet(k)(t-\langle k,\tau\rangle )\big)\,,$$
where $\langle k,\tau\rangle=\suml_{i\in\NN}k_i\tau_i$.

\medskip
The function $\Del$  is well defined since
$$\suml_{k\in\ZZ^\infty}\varphi(k)=\ \ \prodl_{i\in\NN}
\suml_{j\in\ZZ}\mu_i^{|j|}=\ \ \prodl_{i\in\NN}
\frac{1+\mu_i} {1-\mu_i}\ .$$
The latter product is convergent since
$$\suml_{i\in\NN}\mu_i=\suml_{i\in\NN}u^{-1/4}(\tau_i)
<\infty$$
in view of \eqref{eq2.5}.  Since the function
$k\mapsto\varphi(k)\tet(k)$ is bounded on $\ZZ^\infty$  and
since all derivatives $h^{(m)}$ are bounded,
the same argument shows that $\Del$  is a $C^\infty$-function
with
$$\Del^{(m)}(t)=\suml_{k\in\ZZ^\infty}\varphi^{m+1}(k)
\tet^m (k)h^{(m)}\big(\varphi(k)\tet(k)(t-\langle
k,\tau\rangle)\big)\,.$$
Clearly, the function $\Del$  is even.  We have to show
that it satisfies conditions
\eqref{eq2.2}--\eqref{eq2.4}.

\medskip\par\noindent
{\bf Convergence of the integral (\ref{eq2.2}):} Since the
function $h$ is integrable, it suffices to check convergence of
the series
$$
\suml_{k\in\ZZ^\infty}\frac{1}{\tet(k)}=\prodl_{i\in\NN} \
\suml_{\ell\in\ZZ}\frac{1}{\gam^{|\ell|}_{i,\ell}}
\le\prodl_{i\in\NN}\left\{1+2\suml_{\ell\ge
1}|\log\mu_i|^{-\ell}+2\suml_{\ell\ge
1}\mu_i^{\sqrt{\ell}}\right\}.
$$
Since
$$
\suml_{\ell\ge 1}\mu^{\sqrt{\ell}}<\mu+\intl_1^{+\infty}
\mu^{\sqrt{\ell}}d\ell=O\left(\frac{1}{|\log \mu|^2}\right)\,,
$$
for $\mu\to 0$, the right-hand side of the previous expression is bounded
by $\prodl_{i\in\NN}\left\{ 1+\const |\log\mu_i|^{-1}\right\}$. But
this product is finite in view of \eqref{eq2.5}.
\hfill $\qed$

\medskip\par\noindent{\bf Proof of (\ref{eq2.3}):}
Denote by $e^i\in\ZZ^\infty$ the vector with
$e^i_n=\del_{in}$, where $i,n\in\NN$. We have
\begin{eqnarray*}
\Del(t+\tau_i)&=&\suml_{k\in\ZZ^\infty}\varphi(k)h
\big(\varphi(k)\tet(k)(t-\langle k-e^i,\tau\rangle )\big)\\
&=&\suml_{k\in\ZZ^\infty}\varphi(k+e^i)h\big(\varphi(k+e^i)
\tet(k+e^i)(t-\langle k,\tau\rangle)\big)\ .
\end{eqnarray*}
Comparing this with the definition of $\Delta (t)$,
we get that
%(2.8)
\begin{equation}\label{eq2.8}
\frac{\Del(t+\tau_i)}{\Del(t)}\le\sup
\limits_{k\in\ZZ^\infty}\
\frac{\varphi(k+e^i)}{\varphi(k)} \cdot \sup_{k\in\ZZ^\infty}
\sup_{s\in\RR}\frac{h(s)}{h(c_ks)}\ ,
\end{equation}
where
$$c_k=\frac{\varphi(k)\tet(k)}{\varphi(k+e^i)\tet(k+e^i)}\,.$$

\setcounter{theorem}{6}
\begin{lemma}\label{lem2.9}
$$
\mu_i\le\frac{\varphi(k+e^i)}{\varphi(k)}\le
\frac{1}{\mu_i} \qquad \text{and} \qquad
\mu_i\le\frac{\tet(k+e^i)}{\tet(k)}\le\frac{1}{\mu_i}
$$
for all $k\in\ZZ^\infty$, $i\in\NN$.
\end{lemma}

Assume the lemma and note that
$$
\sup_{s\in\RR}\frac{h(s)}{h(cs)}\le 1 \qquad \text{if}\
0 <c\le 1\,,
$$
and
$$
\sup_{s\in\RR}\frac{h(s)}{h(cs)}\le \const\ c(1+
\log^2c) \qquad \text{if}\ c>1\,.
$$
It follows from Lemma~\ref{lem2.9} that $c_k\le \mu_i^{-2}$, and
hence
$$
\sup_{s\in\RR}\frac{h(s)}{h(c_ks)}\le\const\cdot
\left(\frac{1}{\mu^2_i}\right)
\left(1+4\log^2\frac{1}{\mu_i}\right)\,.
$$
Since $\mu_i^{-1}= u(\tau_i)^{1/4}\ge
u(\tau_1)^{1/4}$, we conclude that
$$
\sup_{s\in\RR}
\frac{h(s)}{h(c_ks)}\le \frac{1}{\mu^3_i}
$$
provided $\tau_1$  is sufficiently large. Applying again
Lemma~\ref{lem2.9} and substituting the last inequality into
\eqref{eq2.8} we conclude that
$$\frac{\Del(t+\tau_i)}{\Del(t)}\le \frac{1}{\mu_i}\cdot
\frac{1}{\mu^3_i}=u(\tau_i)\ ,$$
which proves \eqref{eq2.3} modulo the lemma.

\medskip\par\noindent
{\bf Proof of Lemma \ref{lem2.9}:}
The first inequality follows from the fact that
$$\frac{\varphi(k+e^i)}{\varphi(k)}=\mu_i^{|k_i+1|-|k_i|}
=\mu_i^{\pm 1}\ .$$
For the second one, we put
$\alp_{i,\ell}=|\log\mu_i|^{|\ell|}$ and
$\bet_{i,\ell}=\mu_i^{-\sqrt{|\ell|}}$ and notice that
%(2.10)
\setcounter{equation}{7}
\begin{equation}\label{eq2.10}
\frac{\tet (k+e^i)}{\tet(k)}=\frac{\min
(\alp_{i,k_i+1},\bet_{i,k_i+1})}{\min(\alp_{i,k_i},
\bet_{i,k_i})}\ .
\end{equation}
Further,
%(2.11)
\begin{eqnarray}\label{eq2.11}
&&\frac{\alp_{i,\ell +1}}{\alp_{i,\ell}}=
|\log\mu_i|^{\pm 1}\in\left[\mu_i,\frac{1}{\mu_i}
\right]\ ,\\
%(2.12)
\label{eq2.12}
&&\frac{\bet_{i,\ell +1}}{\bet_{i,\ell}}\le\sup_{s\ge 0}
\left(\frac{1}{\mu_i}\right)^{\sqrt{s+1}-\sqrt{s}}=
\frac{1}{\mu_i}\ ,
\end{eqnarray}
and
%(2.13)
\begin{equation}\label{eq2.13}
\frac{\bet_{i,\ell +1}}{\bet_{i,\ell}}
\ge\inf_{s\ge 0}\mu_i^{\sqrt{s+1}-\sqrt{s}}=\mu_i\ .
\end{equation}
Note now that for every 4 positive numbers $a,b,c,d$
$$
\min\left(\frac{a}{c},\frac{b}{d}\right)\ \le\frac{\min
(a,b)}{\min(c,d)}\le\max\left(\frac{a}{c},\frac{b}{d}
\right)\,.
$$
Applying this to \eqref{eq2.10} and using
\eqref{eq2.11}--\eqref{eq2.13}, we conclude that
$$
\mu_i\le\frac{\tet(k+e^i)}{\tet
(k)}\le\frac{1}{\mu_i}\,.
$$
This proves the lemma.\hfill
$\qed$

%\section3
%\section{Verifying asymptotic regularity}\label{sec:veri}
%\setcounter{theorem}{0}

%In this section, we
\medskip
It remains to check that the function $\Del$
satisfies the asymptotic regularity condition (\ref{eq2.4}).
We start with

\medskip\par\noindent{\bf Preliminary estimates: } the function
$\Delta (t)$ satisfies conditions
\begin{equation}
\label{lem3.1} \sup\limits_{|s|\le 1}\ \sup\limits_{t\in\RR}\
\frac{\Del(t+s)}{\Del(t)} <\infty\,,
\end{equation}
and
\begin{equation}
\label{lem3.2} \lim_{t\to\infty} \Del(t) = 0\,.
\end{equation}
Estimate (\ref{lem3.1}) holds for $h$,  and therefore for $\Del$
since the function $k\to\varphi (k)\tet(k)$  is bounded on
$\ZZ^\infty$. Then (\ref{lem3.2}) follows from
integrability of $\Delta$ (see (\ref{eq2.2})\,), and
(\ref{lem3.1}). \hfill $\qed$

The next lemma shows that we did not loose much in the asymptotic
regularity of $\Delta$ compared with that of $h$.

\setcounter{theorem}{13}
\begin{lemma}\label{lem3.5}
For every $m\in\NN$ and every $c\in [0;1) $,
$$
\lim_{t\to\infty}\, \frac{\max_{[t;t+1]}
|\Del^{(m)}|}{\Del^{m+c}(t)} = 0\,.
$$
\end{lemma}

\medskip\par\noindent{\bf Proof:} We show that for every $m\in\NN$
and every $c\in [0;1)$ \setcounter{equation}{14}
\begin{equation}
\label{lem3.3} \lim_{t\to\infty} \frac{\Del^{(m)}}{\Del^{m+c}}(t)
= 0\,.
\end{equation}
Together with (\ref{lem3.1}) this yields the lemma.

In view of (\ref{lem3.2}) it suffices to show that the function
$\Del^{(m)}/\Del^{m+c}$ is bounded on $\RR$  for every $m\in\NN$
and $c\in [0;1)$.  Fix such $m$ and $c$. It is easy to see by
induction in $m$ that
$$
\left( \frac1{t\log^2 t} \right)^{(m)} =\frac{(-1)^m m!}{t^{m+1}
\log^2 t} + \Lambda_m(t)\,, \qquad m\ge 0,
$$
where $\Lambda_m$ is a linear combination of the functions
$\frac1{t^{m+1}\log^k t}$ with $3\le k\le m+2$. Therefore,
$$
h^{(m)}(t) = \frac{(-1)^m m! (1+o(1))}{t^{m+1} \log^2 t}\,, \qquad
t\to +\infty\,,
$$
for each $m\ge 0$, and the function $t\mapsto
h^{(m)}(t)/h^{m+c}(t)$ is bounded. Then we have
$$|\Del^{(m)}(t)|\le\kap\suml_{k\in\ZZ^\infty}\
\varphi^{m+1}(k)\tet^m(k)h^{m+c}(s_k)\,,$$
where $s_k=\varphi(k)\tet(k)(t-\langle k,\tau\rangle )$,
and
$$\kap
=\kap_{m,c}=\sup\limits_{t\in\RR}\frac{|h^{(m)}(t)|}
{h^{m+c}(t)}\,.$$
We claim that
\begin{equation}\label{eq3.4}
\nu_{m,c}:=\sup\limits_{k\in\ZZ^\infty}\
\varphi^{1-c}(k)\tet^m(k)<\infty\ .
\end{equation}
Combining the claim with the elementary inequality
$$\suml_ix^r_i\le \left(\suml_ix_i \right)^r\ ,\ r\ge 1\ ,\ 0\le
x_i\le 1\ ,$$ we readily complete the proof of (\ref{lem3.1}):
\begin{eqnarray*}
&&|\Del^{(m)}(t)|\le\kap\,\nu_{m,c}
\suml_{k\in\ZZ^\infty}\ \varphi^{m+c}(k)h^{m+c}(s_k)\\
&&\le\kap\,\nu_{m,c}\left(\sum_{k\in\ZZ^\infty}\
\varphi(k)h(s_k)\right)^{m+c}=\kap\nu_{m,c}
\Del(t)^{m+c}\ .
\end{eqnarray*}
To prove \eqref{eq3.4}, we set
$$K(m,c)=\{k\in\ZZ^\infty:\varphi^{1-c}(k)\tet^m(k)\le
1\}$$
and check that the complement $\ZZ^\infty\bks K(m,c)$ is
a finite set.  Indeed, if $k\in\ZZ^\infty\bks K(m,c)$,
then
$$
\prodl_{i\in \NN} \left( \mu_i ^{1-c} \gamma_{i,k_i}^m
\right)^{|k_i|} = \varphi^{1-c}(k) \tet^m (k) >1\,,
$$
and therefore at least one of the factors on the left-hand side is
bigger than one. Hence there exists $i\in\NN$ such that
$$\mu^{1-c}_i\gam^m_{i,k_i}\ge 1\ ,$$
which is equivalent to two inequalities:
$$\mu_i^{1-c}|\log\mu_i|^m\ge 1\qquad \text{and}\qquad
\mu_i^{1-c-m/\sqrt{|k_i|}}\ge 1\ .$$
The first inequality shows that
$$
\frac{|\log\mu_i|}{\log|\log\mu_i|}\le\frac{m}{1-c}\,,
$$
therefore there exists a number $j(m,c)$  such that
$i\le j(m,c)$.  The second inequality tells us that
$$
|k_i|\le\left(\frac{m}{1-c}\right)^2\,.
$$
Hence
$$
\#\big(\ZZ^\infty\bks
K(m,c)\big)\le\left(2\left(\frac{m}{1-c}\right)^2+1\right)j(m,c)\,,
$$
and \eqref{eq3.4} follows.  The lemma is proved.\hfill $\qed$

\medskip\par\noindent{\bf Verification of condition (\ref{eq2.4})}:
For a function $v:\RR\to\RR$ denote $(\omega v)(t)=v(t+1)-v(t)$.
Recall that we are proving Property \eqref{eq2.4} which deals with
functions $g_m$ where
$$
g_0(t)=\frac{\Del(t+1)}{\Del(t)}=
\frac{(\omega\Del)(t)}{\Del(t)}+1 \qquad \text{and} \qquad
g_{m+1}(t)= \frac{g'_m(t)}{\Del(t)}\,.
$$
First, note that
$$
g_0(t)-1 = \frac{\Delta(t+1)-\Delta(t)}{\Delta (t)} =
\frac{\Delta'(x_t)}{\Delta (t)}
$$
for some $x_t\in [t,t+1]$. Then Lemma~\ref{lem3.5} yields
$g_0(t)\to 1$ as $t\to\infty$.
It remains to show that $g_m(t)\to 0$  as $t\to\pm\infty$ for
every $m\ge 1$.

\setcounter{theorem}{16}
\begin{lemma}\label{lem3.7}
The function $g_m$  is a finite linear combination of
functions of the form
$$
R=\frac{\omega\Del^{(p)}(\Del')^{\ell_1}\cdots
(\Del^{(m-1)})^{\ell_{m-1}}}{\Del^\ell}\,,
$$
where $p,\ell_1 \nek\ell_{m-1}\ge 0$ and
\setcounter{equation}{18}
$$
%\label{eq3.8}
2\ell_1+\cdots +m\ell_{m-1}+p+2>\ell\ .
\leqno{(3.18)_{\rm m}}
$$
\end{lemma}

\medskip\par\noindent{\bf Proof:} We use induction in $m$.
For $m=0$  we have $p=0$, $\ell=1$, $\ell_1=\cdots =\ell_{m-1}=0$.
Inequality (3.18)$_0$ reads $2>1$. Assume the statement of the
lemma for $m$, and prove it for $m+1$.  Note that $g_{m+1}$ is a
finite linear combination of functions of the form $R'/\Del$. In
turn, $R'/\Del$ is a linear combination of the following
expressions:
$$
\frac{\omega\Del^{(p+1)}(\Del')^{\ell_1}\cdots
(\Del^{(m-1)})^{\ell_{m-1}}}{\Del^{\ell+1}}\,,
$$
$$
\frac{\omega\Del^{(p)}(\Del')^{\ell_1}\cdots
(\Del^{(i-1)})^{\ell_{i-1}}(\Del^{(i)})^{\ell_i-1}
(\Del^{(i+1)})^{\ell_{i+1}+1}(\Del^{(i+2)})^{\ell_{i+2}}
\cdots(\Del^{(m-1)})^{\ell_{m-1}}}{\Del^{\ell +1}}
\,,
$$
where $i=1\nek m-1$, and
$$
\frac{\omega\Del^{(p)}(\Del')^{\ell_1+1}
(\Del^{(2)})^{\ell_2}\cdots (\Del^{(m-1)})^{\ell_{m-1}}}
{\Del^{\ell +2}}\,.
$$
Let us check (3.18)$_{m+1}$ in each of these 3 cases using
(3.18)$_{m}$:
\begin{itemize}
\item[$\bullet$] $\ell +1 <p+1+2+2\ell_1+\cdots +m\ell_{m-1}$;
\item[$\bullet$] $\ell +1 <p+2+\ell_1+\cdots +i\ell_{i-1}+(i+1)
(\ell_i-1)$
$$+(i+2)(\ell_{i+1}+1)+(i+3)\ell_{i+2}+\cdots +
m\ell_{m-1}$$
\item[$\bullet$] $\ell +2 <p+2+2(\ell_1+1)+3\ell_2+\cdots
+m\ell_{m-1}$.
\end{itemize}
This completes the proof. \hfill $\qed$

\medskip\par\noindent
Now we are ready to finish the proof of \eqref{eq2.4}.
It suffices to show that $R(t)\to 0$  as $t\to\infty$,
where $R$ is defined in Lemma~\ref{lem3.7}.
Write
$$
2\ell_1+\cdots +m\ell_{m-1}+p+2=\ell +1+r
$$
with $r\ge 0$.  Choose numbers $\del_0\nek\del_{m-1}$  so
that $\del_i\in [0;1)$ and
$$
(1-\del_0)+\ell_1(1-\del_1)+\cdots
+\ell_{m-1}(1-\del_{m-1})=1\,.
$$
Then
$$\ell =(p+1+\del_0)+\ell_1(1+\del_1)+\ell_2(2+\del_2)
+\cdots +\ell_{m-1}(m-1+\del_{m-1})-r\,.$$
Rewrite $R$  as follows:
$$
R=\frac{\omega\Del^{(p)}}{\Del^{p+1+\del_0}}\cdot
\left(\frac{\Del'}{\Del^{1+\del_1}}\right)^{\ell_1}\cdot\ \cdots\
\cdot
\left(\frac{\Del^{(m-1)}}{\Del^{m-1+\del_{m-1}}}
\right)^{\ell_{m-1}}\cdot\Del^r\,.
$$
Then by Lemma~\ref{lem3.5} and by (\ref{lem3.2}), $R(t)\to 0$ as
$t\to\infty$. This completes the proof of \eqref{eq2.4}, and
therefore finishes off the proof of Theorem \ref{theo2.1}.\hfill
$\qed$

\subsection*{Acknowledgments}
F\"edor Nazarov  generously helped us with the first version of
the growth lemma. In the present form, it appeared after
discussions with Lennart Carleson and Alexei Poltoratskii. Amir
Hadadi (an undergraduate student at Tel Aviv) showed us a
piecewise constant function $\Del$ satisfying \eqref{eq2.2} and a
version of \eqref{eq2.3}. Jean-Christophe Yoccoz sketched an
example of a diffeomorphism of finite smoothness whose growth is
slower than the linear one along a suitable subsequence.  We thank
all of them, as well as Alexander Borichev, Pierre de la Harpe, Anatole Katok,
Patrice Le Calvez and Felix Schlenk for useful discussions and
the referee for helpful critical remarks.

Part of this work was done during the stay of the second named author at
the Mittag-Leffler Institute of the Royal Swedish Academy of Sciences
in the Winter 2002. He thanks this institute for the kind hospitality.

\end{document}